\documentclass[12pt]{amsart}
\usepackage{amscd,amsfonts,amsthm,amssymb,latexsym,amsmath, amscd, amsmath,amsxtra}
\usepackage[all]{xy}
\usepackage{anysize}
\usepackage[sc]{mathpazo} 
\usepackage{soul,xcolor}
\usepackage[toc,page]{appendix}
\usepackage[francais,english]{babel}
\usepackage[utf8]{inputenc}   
\usepackage{hyperref}
\usepackage{stmaryrd}
\usepackage{fancyhdr}
\usepackage{url}
\usepackage{dsfont}

\marginsize{2.5cm}{2.5cm}{2.5cm}{2.5cm}
\newtheorem{theorem}{Theorem}[section] 

\theoremstyle{definition}

\theoremstyle{remark}
\newtheorem{remark}{Remark}

\newcommand{\GL}{\mathrm{GL}}
\newcommand{\SL}{\mathrm{SL}}

\newcommand{\A}{\mathbb{A}}
\newcommand{\La}{\mathrm{L}}

\newcommand{\C}{\mathbb{C}}

\newcommand{\G}{\mathrm{G}}
\newcommand{\U}{\mathrm{U}}
\newcommand{\SU}{\mathrm{SU}}

\begin{document}
\title[Distinction and Base Change]
{Distinction and Base Change}
\author{U. K. Anandavardhanan}

\address{Department of Mathematics, Indian Institute of Technology Bombay, Mumbai - 400076, India.}
\email{anand@math.iitb.ac.in}


\date{}

\begin{abstract}
An irreducible smooth representation of a $p$-adic group $G$ is said to be distinguished with respect to a subgroup $H$ if it admits a non-trivial $H$-invariant linear form. When $H$ is the fixed group of an involution on $G$ it is suggested by the works of Herv\'e Jacquet from the nineties that distinction can be characterized in terms of the principle of functoriality. If the involution is the Galois involution then a recent conjecture of Dipendra Prasad predicts a formula for the dimension of the space of invariant linear forms which once again involves base change. We will describe the proof of this conjecture (in the generic case) for $\SL(n)$ which is joint work with Dipendra Prasad. Then we describe one more newly discovered connection between distinction and base change which is that base change information appears in the constant of proportionality between two natural invariant linear forms on a distinguished representation. This latter result is for discrete series for $\GL(n)$ and is joint with Nadir Matringe. This paper is a report on the author's talk in the International Colloquium on Arithmetic Geometry held in January 2020 at TIFR Mumbai. 
\end{abstract}

\maketitle

\section{Introduction}\label{intro}

The so called relative Langlands program has got a major boost in the last decade thanks to the fundamental works of Gan-Gross-Prasad \cite{ggp12} and Sakellaridis-Venkatesh \cite{sv17}.  Distinguished representations are the basic objects of study in the relative Langlands program and they are studied both locally and globally. In the local setting, given a $p$-adic group $G$, a subgroup $H$ of $G$, and a character $\chi$ of $H$, an irreducible admissible representation $\pi$ of $G$ is said to be $(H,\chi)$-distinguished if
\[{\rm Hom}_H(\pi,\chi)\neq \{0\}.\]
In the global setting, where $H < G$ are adelic groups, and $\pi$ is a cuspidal representation of $G$, the interest is not in any non-trivial $(H,\chi)$-invariant linear form, instead distinction is defined in terms of the non-vanishing of a specific $(H,\chi)$-invariant linear form which is the $(H,\chi)$-period integral. When $\chi$ is the trivial character of $H$, it is often omitted from the definition and we say that the representation is $H$-distinguished. 

In situations such as in \cite{ggp12}, local distinction is connected to the sign of certain local root numbers and global distinction is connected to the special values of certain automorphic $L$-functions. In situations such as in \cite{sv17}, it is expected that distinction can often be characterized in terms of the principle of functoriality.

The notion of a distinguished representation, in both its local and global avatars, goes back to the pioneering work of Harder-Langlands-Rapoport about the Tate conjecture on algebraic cycles in the context of Hilbert-Blumenthal surfaces \cite{hlr86}. In this work they are led to study distinction for the pair $({\rm R}_{F/\mathbb Q}\GL(2),\GL(2))$, where $F$ is a real quadratic number field. Several of the connections that distinction has with other objects of interest are already there in \cite{hlr86}. In particular, it is proved that a cuspidal representation of $\GL_2(\A_F)$ of trivial central character is distinguished precisely when its Asai $L$-function has a pole at $s=1$ and such a representation is characterized as a base change lift of a cuspidal representation of $\GL_2(\A_\mathbb Q)$ with non-trivial central character. Roughly around the same time, the work of Jacquet-Lai introduced the relative trace formula to investigate distinction \cite{jl85}. 

In the nineties, a series of papers of Jacquet and his collaborators investigated distinction for a number of symmetric pairs such as 
\[(\GL(m+n),\GL(m)\times \GL(n)), (\GL(2n),{\rm Sp}(2n)) ~\&~ ({\rm R}_{E/F}{\rm U}(n,E/F),{\rm U}(n,E/F))\]
where $E/F$ is a quadratic extension of $p$-adic fields (or number fields) and ${\rm U}(n,E/F)$ denotes a unitary group defined with respect to $E/F$ . It was proposed that distinction for a symmetric pair $(G,H)$ should have a simple characterization in terms of the principle of functoriality (for instance, see \cite{jy96}). 

An instance of this proposal is that an   irreducible admissible generic representation of $\GL_n(E)$ is distinguished with respect to ${\rm U}(n,E/F)$ if and only if it is a base change lift from $\GL_n(F)$. This is now known by the works of Jacquet \cite{jac04, jac05} and Feigon-Lapid-Offen \cite{flo12} (and \cite{ac89}). Similarly, an   irreducible admissible generic representation of $\GL_n(E)$ is $(\GL_n(F),\omega_{E/F}^{n-1})$-distinguished if and only if it is a base change lift from ${\rm U}(n,E/F)$, where $\omega_{E/F}$ is the quadratic character of $F^\times$ associated to $E/F$. This suggestion is sometimes referred to as the Flicker-Rallis conjecture (it is the local version of the global conjecture stated in \cite[p. 143]{fli91}). This too is now known by combining \cite[Theorem 5.2]{mat11} and \cite[Lemma 2.2.1]{mok15}.    

This paper is an informal exposition of two results of the author which bring to light the close relationship between distinction and base change; one of these is joint work with Dipendra Prasad \cite[Theorem 5.6]{ap18} and the other is joint work with Nadir Matringe \cite[Theorem 6.1]{am20}. The first is Theorem \ref{ap} and the second is Theorem \ref{am} of this paper. Theorem \ref{ap} proves a conjecture of Dipendra Prasad \cite[Conjecture 13.3]{pra20} for $\G=\SL(n)$. We start by recalling a few aspects of Prasad's conjecture in Section \ref{prasad} and then summarize a number of facts connecting distinction for $(\GL_n(E),\GL_n(F))$ and base change from $\U(n,E/F)$ to $\GL_n(E)$ (cf. \cite{mok15}) in Section \ref{mok}. In Section \ref{sln}, we take up Theorem \ref{ap}, and we discuss Theorem \ref{am} in Section \ref{fd}.

\section{Prasad's Conjecture}\label{prasad}

Prasad's conjecture is for $p$-adic fields and for a Galois pair \cite[Conjecture 13.3]{pra20}. As in Section \ref{intro}, $E/F$ is a quadratic extension of $p$-adic fields. Let $\G$ be a connected reductive group defined over $F$. Let $G=\G(E)$ and $H=\G(F)$. The Galois involution $\sigma$ acts on $G$ and hence on representations of $G$. For an irreducible admissible representation $\pi$ of $G$ its Galois conjugate is denoted by $\pi^\sigma$. The contragredient representation of $\pi$ is denoted by $\pi^\vee$. If $\pi$ is such that $\pi^\vee \cong \pi^\sigma$, we say that $\pi$ is conjugate self-dual.

There are two objects associated to the data $(\G,E/F)$ which appear in the formulation of the conjecture. One of these is the opposition group denoted by $\G^{\rm op}$ constructed in \cite[\S 5]{pra20} which in particular has the property that it is isomorphic to $\G$ over $E$. The other is a character of $H$ of order $\leq 2$, denoted by $\omega_\G$,  constructed in \cite[\S 6]{pra20}. 

The conjecture has three parts. The first part asserts that an irreducible admissible representation of $G$ which is $(H,\omega_G)$-distinguished is conjugate self-dual and moreover its $\La$-packet arises as the base change lift of an $\La$-packet of $\G^{\rm op}(F)$.  

The two examples considered in Section \ref{intro} correspond to 
\begin{itemize}
\item $\G=\GL(n), \G^{\rm op} = {\rm U}(n,E/F), \omega_G = \omega_{E/F}^{n-1}\circ \det$,
\item $\G={\rm U}(n,E/F), \G^{\rm op}=\GL(n), \omega_G=1$.
\end{itemize}
Thus, we have already seen that the first part of Prasad's conjecture is true in the case of both these examples for generic representations. In fact the assertion on conjugate self-duality is known to be true for all irreducible admissible representations. For $\G=\GL(n)$, this is \cite[Proposition 12]{fli91}, and for $\G={\rm U}(n,E/F)$, this is \cite[Theorem 0.2]{flo12}.   

The second and third parts of the conjecture are for generic representations so assume that $\G$ is quasi-split over $F$. Fix a Borel subgroup $B$ of $\G$ and let $N$ be its unipotent radical.

The second part of the conjecture probes for distinction inside a generic $\La$-packet and proposes a simple recipe to detect distinction. It says that an irreducible admissible representation which is generic for a non-degenerate character of $N(E)/N(F)$ is $(H,\omega_\G)$-distinguished provided its $\La$-packet is a base change lift of an $\La$-packet of $\G^{\rm op}(F)$.    

The third and perhaps the most important part of the conjecture is the multiplicity formula. We do not state it as it involves more technical terms and refer to \cite[Conjecture 13.3]{pra20} for the precise statement. It suffices to say that a key ingredient in the conjectural formula for the multiplicity - which is $\dim_{\mathbb C} {\rm Hom}_H(\pi,\omega_G)$ - is the cardinality of the fiber of the base change map from $\La$-packets of $\G^{\rm op}(F)$ to $\La$-packets of $G$. Thus the connection between distinction and base change is indeed quite deep. 

For $\G=\GL(n)$, the corresponding Galois pair is a Gelfand pair \cite[Proposition 11]{fli91}; i.e.,
\[\dim_{\mathbb C} {\rm Hom}_H(\pi,\omega_G) = 1\]
for every irreducible admissible representation $\pi$ of $G$ which is $(H,\omega_G)$-distinguished. This fits well with the conjectural multiplicity formula which in this case equals the cardinality of the fiber of the base change map from ${\rm U}(n,E/F)$ to $\GL_n(E)$ as this map is known to be injective (cf. \cite[Proposition 7.10]{pra20}).  

For $\G={\rm U}(n,E/F)$, the multiplicity can be more than $1$. For instance, it is easy to see that the principal series representation Ps$(\chi,\chi^{-1})$ of $\GL_2(E)$ with $\chi=\chi^\sigma$ is both $\GL_2(F)$-distinguished and $(\GL_2(F),\omega_{E/F})$-distinguished and that the corresponding invariant linear functionals are ${\rm U}(2,E/F)$-invariant. It can be seen that 
\[\dim_{\mathbb C} {\rm Hom}_H(\pi,\omega_G) = 2\]
in this case. More generally, if $\pi$ is an   irreducible admissible generic representation of $\GL_n(E)$ then it is parabolically induced from a number of essentially square-integrable representations $\pi_i$ of $\GL_{n_i}(E)$, $1 \leq i \leq t$, $n=n_1+\dots+n_t$, say
\[\pi = \pi_1 \times \ldots \times \pi_t.\]
Suppose $\pi$ and $r$ many of these $\pi_i$'s are Galois invariant. The conjectural multiplicity formula would then predict
\[\dim_{\mathbb C} {\rm Hom}_H(\pi,\omega_G) = 
\begin{cases}
2^{r-1} &\text{$r\geq 1$,} \\
1 &\text{$r=0$.}
\end{cases}
\]
We refer to \cite[\S 17]{pra20} for more details. By \cite[Theorem 0.2]{flo12}, the right hand side is known to be a lower bound with equality known to hold if the Galois invariant $\pi_i$'s are distinct. Recently, R. Beuzart-Plessis has proved this multiplicity formula in general \cite[Theorem 3]{bp20}.

\section{Base change from ${\rm U}(n,E/F)$ to $\GL_n(E)$}\label{mok}

In Section \ref{prasad}, we have briefly mentioned the connection between distinction for the Galois pair $(\GL_n(E),\GL_n(F))$ and base change from a unitary group. We make it more precise in this section. We closely follow \cite[\S 2.1, 2.2]{mok15}.

In order to discuss base change we need to introduce the notion of a Langlands parameter of a group $\G$ defined over a $p$-adic field $k$. It is an admissible homomorphism from the Weil-Deligne group of $k$ 
\[W_k^\prime = W_k \times \SL_2(\mathbb C)\]
to the Langlands dual group 
\[{^L}\G = \G^\vee \rtimes W_k\]
of $\G$, where $\G^\vee$ is the complex dual group of $\G$. Thus, a Langlands parameter is a continuous homomorphism
\[\varphi: W_k^\prime \rightarrow {^L}\G\]
that commutes with the natural projections $W_k^\prime \rightarrow W_k$, ${^L}\G \rightarrow W_k$, it sends $W_k$ to semisimple elements and its restriction to $\SL_2(\mathbb C)$ is algebraic. Two Langlands parameters are equivalent if they are conjugate by $\G^\vee$. We denote by $\Phi(\G)$ (or sometimes by $\Phi(\G(k))$) the set of equivalence classes of Langlands parameters of $\G$.     

The groups of interest to us are $\GL(n)$ and the quasi-split unitary group defined by  
\[\U(n) = \{g \in \GL_n(E) \mid {^t}g^\sigma J g = J \},\] 
where 
\[J = \mbox{anti-diag}((-1)^{n-1},\dots,1)) \]
is the anti-diagonal matrix with $\pm 1$ alternating.  The Langlands dual group of $\GL(n)$ is $\GL_n(\mathbb C)$. The Langlands dual group of $\U(n)$ is given by   
\[{^L}\U(n) = \GL_n(\mathbb C) \rtimes W_F,\] 
where the Weil group $W_F$ of $F$ acts by projection to Gal$(E/F)$,  and $w_\sigma \in W_F \smallsetminus W_E$ acts as the automorphism
\[g \mapsto J~ {^t}g^{-1}J^{-1}.\]

The base change map from $\U(n)$ to $\GL_n(E)$ at the level of Langlands parameters is a map from 
\[\Phi(\U(n)) = \{\varphi:W_F^\prime \rightarrow {^L}\U(n)\}/\sim \] to 
\[\Phi(\GL_n(E)) = \{\rho:W_E^\prime \rightarrow \GL_n(\mathbb C)\}/\sim \] 
given by
\begin{center}
\begin{tabular}{rcccc}
BC & : & $\Phi(\U(n))$ & $\rightarrow$ & $\Phi(\GL_n(E))$ \\ 
 & & $\varphi$ & $\mapsto$ & $\varphi |_{_{W_E^\prime}}.$
\end{tabular}
\end{center} 

Now the Flicker-Rallis conjecture (see \cite[p. 143]{fli91} where it is stated in the global context) is the assertion that, for an   irreducible admissible generic representation $\pi$ of $\GL_n(E)$ with Langlands parameter $\rho_\pi$,  
\begin{align*}
\rho_\pi \in \mbox{Image(BC)} \iff 
\begin{cases}
\pi \mbox{ is $\GL_n(F)$-distinguished} &\text{for $n$ odd,} \\ 
\pi \mbox{ is $(\GL_n(F),\omega_{_{E/F}})$-distinguished} &\text{for $n$ even.}  
\end{cases} 
\end{align*}

\begin{remark}\label{unstable}
What we described above is what is called the stable base change map from $\U(n)$ to $\GL_n(E)$. There is also an unstable base change map with respect to an extension $\kappa$ of $\omega_{E/F}$ to $E^\times$. A representation $\pi$ of $\GL_n(E)$ is in the image of the stable base change map if and only if $\pi \otimes \kappa$ is in the image of the unstable base change map with respect to $\kappa$.
\end{remark}

As mentioned in Section \ref{prasad}, this is now proved thanks to \cite[Theorem 5.2]{mat11} and \cite[Lemma 2.2.1]{mok15}. We end this section by stating these two results and by paraphrasing the Flicker-Rallis conjecture in the language of Prasad's conjecture.

To this end, we introduce the notion of parity of a conjugate self-dual Langlands parameter for $\GL_n(E)$. So let 
\[\rho:W_E^\prime \rightarrow \GL_n(\mathbb C)\]
be such that $\rho^\sigma \cong \rho^\vee$, where $\rho^\sigma(g) = \rho(w_\sigma^{-1}g w_\sigma)$ for $w_\sigma \in W_F \smallsetminus W_E$. Then $\rho$ is said to be of parity $\eta(\rho) \in \{\pm 1\}$ if there is a non-degenerate bilinear form $B$ on $\rho$ with 
\begin{itemize}
\item $B(\rho(g)v,\rho^\sigma(g)w) = B(v,w)$, 
\item $B(v,w) = \eta(\rho) \cdot B(w,\rho(w_\sigma^2)v)$.
\end{itemize}

Thus the Flicker-Rallis conjecture follows from combining Theorem \ref{matringe} and Theorem \ref{mok}. 

\begin{theorem}[Matringe]\label{matringe}
An   irreducible admissible generic representation of $\GL_n(E)$ is $\GL_n(F)$-distinguished if and only if its Langlands parameter is conjugate self-dual of parity $+1$.  
\end{theorem} 

Theorem \ref{matringe} is \cite[Theorem 5.2]{mat11} which builds on earlier results in \cite{kab04,akt04,ar05,mat09}. Theorem \ref{mok} is \cite[Lemma 2.2.1]{mok15}.

\begin{theorem}[Mok]\label{mok}
A Langlands parameter for $\GL_n(E)$ is in the image of the base change map from $\U(n)$-parameters if and only if it is conjugate self-dual of parity $(-1)^{n-1}$.
\end{theorem}

Now we state all these together in the language of Prasad's conjecture.

\begin{theorem}[Flicker-Matringe-Mok]\label{fmm}
An   irreducible admissible generic representation $\widetilde{\pi}$ of $\GL_n(E)$ is $(\GL_n(F),\omega_{E/F}^{n-1})$-distinguished if and only if its Langlands parameter $\widetilde{\rho}_{\widetilde{\pi}}$ is in the image of
\[{\rm BC}: \Phi(\U(n)) \rightarrow \Phi(\GL_n(E)),\]
and moreover,
\begin{equation*}
\dim_\C {\rm Hom}_{GL_n(F)}(\widetilde{\pi},\omega_{_{E/F}}^{n-1}) = |{\rm BC}^{-1}(\widetilde{\rho}_{\widetilde{\pi}})|. 
\end{equation*}
\end{theorem}
 
\section{$(\SL_n(E),\SL_n(F))$}\label{sln}

Prasad's conjecture for $\SL(n)$ is proved in \cite{ap18} (see also \cite{ap03,ana05} for some of the early works in this direction). In this section, following \cite{ap18}, we summarise the key steps involved in its proof.

The opposition group for $\G=\SL(n)$ over $E/F$ is $\G^{\rm op}=\SU(n,E/F)$ and Prasad's character $\omega_G$ is the trivial character. Thus, we are interested in the space of ${\SL_n(F)}$-invariant linear forms on an   irreducible admissible generic representation $\pi$ of $\SL_n(E)$ and also in the base change map from $\SU(n)$ to $\SL_n(E)$.

Base change for $\SU(n)$ fits into the commutative diagram   
\[ 
\xymatrix@C=2pc@R=2pc{
\Phi(\SU(n)) \ar[r]^{p{\rm BC}} & \Phi(\SL_n(E)) \\
\Phi(\U(n)) \ar[u]_{p_F} \ar[r]^{{\rm BC}} & \Phi(\GL_n(E)) \ar[u]_{p_E}
}
\]   
where $p_F$ and $p_E$ are the natural projections induced by the homomorphism 
\[\GL_n(\mathbb C) \rightarrow {\rm PGL}_n(\mathbb C).\]
    
The key observation is that the maps $p_F, p_E$ are surjective and BC is injective. The surjectivity of $p_E$ follows from Tate's theorem according to which $H^2(W_E^\prime,\mathbb C^\times)=0$ where $W_E^\prime$ acts trivially on $\mathbb C^\times$. The surjectivity of $p_F$ also follows from a similar second cohomology vanishing result (cf. \cite[Lemma 5.1]{ap18}). The map BC is injective is the result \cite[Proposition 7.10]{pra20}.
  
To proceed, let $\rho \in \Phi(\SU(n))$.    Let $\widetilde{\rho} \in p_F^{-1}(\rho)$. Now observe that
\[{\rm BC}(p_F^{-1}(\rho))    = \{{\rm BC}(\widetilde{\rho}) \otimes \chi \mid \chi \in \widehat{E^\times/F^\times}\},\]   
and
\[p_E^{-1}(p{\rm BC}(\rho))    = \{{\rm BC}(\widetilde{\rho}) \otimes \chi \mid \chi \in \widehat{E^\times}\}.\]

The above observation leads us to make the following crucial definitions (see also \cite[Remark 2]{ana05}). Two members of BC$(\Phi(\U(n)))$ are \emph{weakly (resp. strongly) equivalent} if they differ by a character of $E^\times$    (resp. $E^\times/F^\times$). Strong and weak equivalence classes are similarly defined in the set of $(\GL_n(F),\omega_{E/F}^{n-1})$-distinguished representations.  

It follows that the cardinality of 
\[\{\mu \in \Phi(\SU(n)) \mid p{\rm BC}(\mu) = p{\rm BC}(\rho)\}\] 
is the number of strong equivalence classes in the weak equivalence class of BC$(\widetilde{\rho})$.

We now state the main theorem of \cite{ap18} which is the exact analogue of Theorem \ref{fmm} for $\SL(n)$    with an additional condition to probe for distinction inside an $\La$-packet \cite[Theorem 5.6]{ap18}. 

\begin{theorem}[Anandavardhanan-Prasad]\label{ap}
An   irreducible admissible generic representation $\pi$ 
of $\SL_n(E)$ is distinguished by $\SL_n(F)$ if and only if 
\begin{enumerate}
\item[(1)] its Langlands parameter $\rho_\pi$ is in the image of 
\[p{\rm BC}: \Phi(\SU(n)) \rightarrow \Phi(\SL_n(E)),\]
\item[(2)] $\pi$ has a Whittaker model for a non-degenerate character of $N(E)/N(F)$.
\end{enumerate}
Further, if ${\rm Hom}_{\SL_n(F)}(\pi,1) \not = 0$,
\begin{equation*}
\dim_\C {\rm Hom}_{\SL_n(F)}(\pi,1) = |p{\rm BC}^{-1}(\rho_\pi)|.
\end{equation*}
\end{theorem}

Part (1) of Theorem \ref{ap} follows from the commutative diagram above for base change from $\SU(n)$ together with the fact that $\pi \in \widetilde{\pi}|_{_{\SL_n(E)}}$ for some irreducible admissible generic representation $\widetilde{\pi}$ of $\GL_n(E)$ which can be taken to be $(\GL_n(F),\omega_{E/F}^{n-1})$-distinguished.

The key ingredient in proving (2) is Theorem \ref{fo} below, which follows from combining a number of results. Firstly, we have a result due to Flicker by which for an irreducible admissible unitary generic representation $\pi$ of $\GL_n(E)$ and for a non-degenerate character $\psi$ of $N(E)/N(F)$, the integral
\[\int_{N(F)\backslash P(F)} W(p) dp\]
is absolutely convergent for $W$ in the $\psi$-Whittaker model $\mathcal W(\pi,\psi)$ of $\pi$, where $P(F)$ denotes the mirabolic subgroup of $\GL_n(F)$ \cite[Lemma in \S 4]{fli88}. Also, by \cite[Proposition in \S 4]{fli88}, the linear form defined on $\mathcal W(\pi,\psi)$ by such an integral is non-trivial. If $\pi$ is non-unitary but generic and also $\GL_n(F)$-distinguished then it is shown in \cite[Section 7]{am17} that the integral
\[\int_{N(F)\backslash P(F)} W(p)|\det p|^{s-1} dp\]
(which is convergent for Re$(s)$ large and admitting a meromorphic continuation to $\mathbb C$ by the Rankin-Selberg theory of the Asai $L$-function and in particular by \cite[Lemma 3.5]{mat09})   
is holomorphic at $s=1$. We denote this regularized integral by 
\[\int^*_{N(F)\backslash P(F)} W(p) dp\]
which is not identically zero as a linear form on $\mathcal W(\pi,\psi)$ (see, for instance, \cite[Theorem 7.2]{am17}).
Secondly, we have a result due to Youngbin Ok by which 
\[{\rm Hom}_{P(F)}(\pi,1) = {\rm Hom}_{\GL_n(F)}(\pi,1)\]
for any irreducible admissible $\GL_n(F)$-distinguished representation $\pi$
(see \cite[Proposition 2.1]{mat10} and \cite[Theorem 3.1]{off11}). We remark here that Ok's result in the tempered case is proved independently in \cite[Theorem 1.1]{akt04}.
\begin{theorem}\label{fo}
The unique, up to multiplication by scalars, $\GL_n(F)$-invariant linear form on a $\GL_n(F)$-distinguished   irreducible admissible generic representation $\widetilde{\pi}$ of $\GL_n(E)$ is given on its $\psi$-Whittaker model $\mathcal W(\pi,\psi)$ by
\[\ell(W) = \int^*_{N(F)\backslash P(F)} W(p) dp,\]
where $\psi$ is a non-degenerate character of $N(E)/N(F)$.
\end{theorem}

In order to prove the multiplicity formula
\[\dim_\C {\rm Hom}_{\SL_n(F)}(\pi,1) = |p{\rm BC}^{-1}(\rho_\pi)|,\]
choose an   irreducible admissible generic representation $\widetilde{\pi}$ of $\GL_n(E)$ containing the representation $\pi$ which is $(\GL_n(F),\omega_{E/F})$-distinguished and $\widetilde{\rho}$ such that BC$(\widetilde{\rho})=\widetilde{\rho}_{\widetilde{\pi}}$. The right hand side is the number of strong equivalence classes in the weak equivalence class of $\widetilde{\rho}_{\widetilde{\pi}}$    (inside Image(BC)). The left hand side can be shown to be    the number of strong equivalence classes in the weak equivalence class of $\widetilde{\pi}$  (among the $(\GL_n(F),\omega_{E/F}^{n-1})$-distinguished irreducible generic representations of $\GL_n(E)$).

\begin{remark}\label{rem1}
A consequence of Theorem \ref{ap} (2) is that an $\SL_n(F)$-distinguished   irreducible admissible generic representation of $\SL_n(E)$ is conjugate self-dual. This follows from the uniqueness of non-degenerate Whittaker models for $\GL_n(E)$. Indeed, for an irreducible $\psi$-generic representation $\pi$ of $\SL_n(E)$ which is $\SL_n(F)$-distinguished let $\widetilde{\pi}$ be an   irreducible admissible generic representation of $\GL_n(E)$ such that $\pi$ appears in its restriction to $\SL_n(E)$. By multiplicity one of Whittaker functionals on $\widetilde{\pi}$ we see that $\pi$ is the unique $\psi$-generic representation in its $\La$-packet. A character $\psi$ of $N(E)/N(F)$ has the property that $\psi^{-1}=\psi^\sigma$ and therefore it follows that both $\pi^\vee$ and $\pi^\sigma$ are generic with respect to the same character of $N(E)$ and thus they are isomorphic.
\end{remark}

\begin{remark}\label{rem2}
In \cite{am21}, Theorem \ref{ap} (2) is proved for any unitary representation of $\SL_n(E)$ and therefore distinction by $\SL_n(F)$ implies conjugate self-duality for any unitary representation of $\SL_n(E)$ by the multiplicity one result for degenerate Whittaker models by the argument outlined in Remark \ref{rem1}. Thus we have established one of the assertions in Prasad's conjecture for all the unitary representations of $\SL_n(E)$. The role played by Theorem \ref{fo} in the proof of Theorem \ref{ap} (2) is played by \cite[Propositions 2.2 and 2.5]{mat14} in the unitary context.
\end{remark}

\begin{remark}
The assertions in Remark \ref{rem1} are true also over finite fields; i.e., for the pair $(\SL_n(\mathbb F_{q^2}),\SL_n(\mathbb F_q))$. This is \cite[Theorem 5.1]{am20}. 
\end{remark}

\section{A new connection between distinction and base change}\label{fd}

In Section \ref{prasad}, we saw Prasad's conjecture which related distinction and base change in a precise way and in particular the cardinality of the fiber of the base change from the opposition group $\G^{\rm op}(F)$ goes into the conjectural formula for the dimension of the space of $(\G(F),\omega_\G)$-invariant forms on an   irreducible admissible generic representation of $\G(E)$. In Section \ref{sln}, we saw the main ideas behind the proof of Prasad's conjecture for $\SL(n)$ in \cite{ap18}. 

In this section, we present a result, joint with Nadir Matringe, which illustrates the connection between distinction and base change in yet another way which is that base change information appears in the constant of proportionality between two natural invariant linear forms on a distinguished representation. The result is for the pair $(\GL_n(E),\GL_n(F))$ and for $\GL_n(F)$-distinguished discrete series representations of $\GL_n(E)$ and it is contained in \cite[Section 6]{am20}. We give an informal introduction to the result and its proof.

The main points to keep in mind are: 
\begin{enumerate}
\item The pair $(\GL_n(E),\GL_n(F))$ is of multiplicity one \cite[Proposition 11]{fli91}.  
\item There are two natural $\GL_n(F)$-invariant forms on a $\GL_n(F)$-distinguished discrete series representation of $\GL_n(E)$. One due to Flicker \cite[Lemma in \S4]{fli88}, denoted by $\ell$, which we saw in Theorem \ref{fo}, and the other due to Kable, say $\lambda$ \cite[Theorem 4]{kab04}. 
\item\label{3} The distinguishing linear forms $\lambda$ and $\ell$ differ by a constant by Flicker's multiplicity one result.
\item Flicker-Rallis conjecture, recalled in Theorem \ref{fmm}, according to which distinction for $(\GL_n(E),\GL(_n(F))$ is related to base change from $\U(n,E/F)$. 
\end{enumerate}  

Our result evaluates the proportionality constant in (\ref{3}) above and it involves the formal degrees of
the base changed and base changing representations. We state the result towards the end of this section (cf. Theorem \ref{am}).   

Thus the two inputs for the statement of Theorem \ref{am}  are base change and formal degree.  We have introduced base change from $\U(n)$ to $\GL_n(E)$ in Section \ref{mok}. Let us now recall the definition of the formal degree of a discrete series representation.  
 
If $\pi$ is a discrete series representation of a $p$-adic group, there exists $d_\mu(\pi) \in \mathbb R_{>0}$ such that
\[\int_{G/Z} \langle \pi(g)v, v^\prime \rangle \overline{\langle \pi(g)w, w^\prime \rangle}d\mu(g) = \frac{1}{d_\mu(\pi)} \langle v, w \rangle \overline{\langle v^\prime, w^\prime \rangle} . \]  
This $d_\mu(\pi)$ is the formal degree of $\pi$, which of course depends on the choice of the Haar measure $\mu$.

\begin{remark}
Recall the orthogonality relations for matrix coefficients for a finite group.   A matrix coefficient of a (unitary) representation $\pi$ is a function on $G$ given by
\[g \mapsto \langle \pi(g)v,v^\prime \rangle \]
where $v, v^\prime \in \pi$.   If $\pi$ is irreducible,   we have
\[\frac{1}{|G|} \sum_{g \in G} \langle \pi(g)v, v^\prime \rangle \overline{\langle \pi(g)w, w^\prime \rangle} = \frac{1}{\dim \pi} \langle v, w \rangle \overline{\langle v^\prime, w^\prime \rangle} . \]
Thus, formal degree, for infinite dimensional representations, plays the role of the dimension of a finite dimensional representation.
\end{remark}

The Hiraga-Ichino-Ikeda conjecture gives a formula for the formal degree of a discrete series representation of a $p$-adic group in terms of its adjoint gamma function evaluated at $s=0$ \cite[Conjecture 1.4]{hii08}. This conjecture is proved for $\GL(n)$ in \cite{hii08} itself (cf. \cite[Theorem 3.1]{hii08}) and it is proved for $\U(n,E/F)$ by Beuzart-Plessis \cite[Corollary 5.5.4]{bp18b}.

With respect to the specific choice of Haar measure as in [HII '08], we have:  
\begin{enumerate}
\item\label{hii-1} For a discrete series representation $\pi$ of $\GL_n(E)$, its formal degree is given by \cite[Theorem 3.1]{hii08}
\[d(\pi) =\frac{1-q^{-1}}{n} \left| \displaystyle{\lim_{s \rightarrow 0}} \frac{\gamma(s,\pi \times \pi^\vee, \psi)}{1-q^{-s}} \right|,\]
where the gamma factor is the Rankin-Selberg gamma factor.
\item\label{hii-2} For a discrete series representation $\rho$ of $\U(n,E/F)$, its formal degree is given by \cite[Corollary 5.5.4]{bp18b} (see also \cite[Lemma 8.1]{hii08})
\[d(\rho) = \frac{1}{2} \left| \gamma(0,\pi,r^\prime,\psi_0) \right|,\]
where the gamma factor is the twisted Asai gamma factor.
\end{enumerate}

\begin{remark}
There are three ways of defining a gamma factor. Via the Langlands formalism, via the Langlands-Shahidi method, and via the Rankin-Selberg integral method. All these three definitions coincide in the cases considered in this section. For pairs of representations of $\GL(n)$, this follows from the local Langlands conjecture and by \cite[Theorem 5.1]{sha84}. In case of $\U(n)$, this follows from \cite[Theorem 1]{sha18} and \cite[Theorem 3.4.1]{bp18a} (see also \cite[Remark 3.5]{ar05}). In the latter case, an appropriate normalization is required for the definition via the integral method (cf. \cite[Definition 9.10]{akmss18} and \cite[Theorem 3.4.1]{bp18a}).  
\end{remark}

Let $\psi$ be a non-degenerate character of $N(E)/N(F)$. Consider  
\[\ell(W) = \int_{N(F) \backslash P(F)} W(p) dp\]  
on the $\psi$-Whittaker model $\mathcal W(\pi,\psi)$ of an irreducible unitary generic representation $\pi$ of $\GL_n(E)$ (cf. Theorem \ref{fo}).  As mentioned earlier, this linear form is first considered by Flicker who also proved the absolute convergence of the above integral \cite[\S 4]{fli88}. The form $\ell$ is always non-zero  and clearly $P(F)$-invariant. It is $\GL_n(F)$-invariant precisely when $\pi$ is $\GL_n(F)$-distinguished (see Section \ref{sln}).
 
Now consider
\[\lambda(W) = \int_{F^\times N(F) \backslash \GL_n(F)} W(g) dg\]  
on $\mathcal W(\pi,\psi)$ which is obviously $\GL_n(F)$-invariant but this integral is convergent only when $\pi$ is a discrete series representation \cite[Lemma 2]{kab04}. So assume $\pi$ is a discrete series representation. The form $\lambda$ is non-zero precisely when $\pi$ is $\GL_n(F)$-distinguished \cite[Theorem 4]{kab04}. 

Now \cite[Theorem 6.1]{am20} is the following.

\begin{theorem}[Anandavardhanan-Matringe]\label{am}
Let $\pi$ be a discrete series representation of $\GL_n(E)$ which is $\GL_n(F)$-distinguished. Let $\rho$ be the (discrete series) representation of $\U(n,E/F)$ that base changes to $\pi$ (stably or unstably depending on the parity of $n$). Let $d(\rho)$ (resp. $d(\pi)$) denote the formal degree of $\rho$. Then,
\[\lambda = c \cdot \frac{d(\rho)}{d(\pi)} \cdot \ell,\]
where $c$ is a positive constant that does not depend on the representations $\rho$ and $\pi$.
\end{theorem}

For the proof of Theorem \ref{am}, we refer to \cite[Section 6]{am20}. However, we indicate the key ingredients in its proof in a sort of informal way.

The starting point of the proof of Theorem \ref{am} is the functional equation for the Asai $L$-function defined via the Rankin-Selberg integral method. This is due to \cite[Appendix]{fli93} and \cite[Proposition 2]{kab04}. Thus, if
\[Z(s,W,\Phi)=\int_{N_n(E)\backslash \GL_n(E)} W(g)\Phi((0,\dots,0,1)g)|\det g|^s dg\]
for $W \in \mathcal W(\pi,\psi)$ and for a Schwartz-Bruhat function on $E^n$, we have 
\begin{align}\label{fe}
Z(1-s,\widetilde{W},\widehat{\Phi})=\gamma(s,\pi,r,\psi)Z(s,W,\Phi)
\end{align}
where $\widetilde{W} \in \mathcal W(\pi^\vee,\psi^{-1})$ is given by $W(J {^t}g^{-1})$ and $\widehat{\Phi}$ is the Fourier transform of $\Phi$. 

Now by the proof of \cite[Theorem 4]{kab04}, the right hand side of (\ref{fe}) is connected to the linear form $\lambda$.

On the other hand, by \cite[Theorem 1.4]{akt04}, the left hand side of $(\ref{fe})$ is connected to the linear form $\ell$. This step is subtle and involves a new functional equation \cite[Theorem 6.3]{am17} which in turn requires knowledge of the sign of the local root number which is a particular case of \cite[Theorem 3.6]{mo17} (cf. \cite[Conjecture 5.1]{ana08}). 

Thus, the linear forms $\lambda$ and $\ell$ get connected via $\gamma(0,\pi,r,\psi)$, and via the factorization
\[\gamma(s,\pi \times \pi^\vee, \psi)= \gamma(s,\pi,r,\psi)\gamma(s,\pi,r^\prime,\psi)\]
they get connected to the formal degrees $d(\pi)$ and $d(\rho)$ by the Hiraga-Ichino-Ikeda conjecture for $\GL(n)$ as well as $\U(n)$.  

This analysis finally leads to the relation
\[\lambda = c \cdot \epsilon(1/2,\pi,r,\psi)\cdot \frac{d(\rho)}{d(\pi)} \cdot \ell,\]
where $c$ is a measure theoretic positive constant which does not depend on the representations $\pi$ and $\rho$.

Theorem \ref{am} follows since in the situation at hand it is known that 
\[\epsilon(1/2,\pi,r,\psi)=1\]
by combining the results \cite[Theorem 1.1]{ana08}, \cite[Theorem 1]{sha18}, and \cite[Theorem 3.4.1]{bp18a}.

\begin{remark}
In \cite{am20}, we were led to the formulation of Theorem \ref{am} by a result over finite fields showing that, for $\G=\GL(n)$ or $\G=\U(n,q)$, the Bessel function $B_{\pi,\psi}$ is a test vector for the natural $\G(\mathbb F_q)$-invariant linear form on an irreducible generic representation $\pi$ of $\G(\mathbb F_{q^2})$ which is a base change from $\G^{\rm op}(\mathbb F_q)$. Here, $\psi$ is a non-degenerate character of $N(\mathbb F_{q^2})/N(\mathbb F_q)$, where $N$ is the unipotent radical of a (fixed) Borel subgroup of $\G$. In fact, by \cite[Theorem 1.1]{am20},
\[\frac{1}{\G(\mathbb F_q)} \sum_{h \in \G(\mathbb F_q)} B_{\pi,\psi}(h) = \frac{\left| \G(\mathbb F_{q^2})\right|}{\left| \G(\mathbb F_q)\right| \left| \G^{\rm op}(\mathbb F_q)\right|} \cdot \frac{\dim \rho}{\dim \pi},\]
where $\rho$ is the irreducible generic representation of $\G^{\rm op}(\mathbb F_q)$ that base changes to $\pi$. This identity for finite fields is generalized to all irreducible generic uniform representations for any connected quasi-split reductive group by Chang Yang \cite[Theorem 1.1]{yan20}. It will be of interest to look for a $p$-adic analogue as in Theorem \ref{am} for $\G=\GL(n)$.     
\end{remark}

\section*{Acknowledgements} 

My own results described in this paper are from joint work with Nadir Matringe and Dipendra Prasad. It has been a pleasure to work together and I thank both of them heartily. Also, I am grateful to both of them for their suggestions on the first draft of this manuscript. I would like to thank Dipendra Prasad for constant support and encouragement over the last two decades. Indeed, our work on $\SL(n)$ originated when I was a postdoc at TIFR Mumbai during 2003-2005. I would like to thank C. S. Rajan, my initial work with him, again during my postdoc period, has had a major impact on several of the results presented in this paper. Thanks are due to Rajat Tandon for introducing me to the beautiful area of representation theory of $p$-adic groups and in particular to distinguished representations. Finally I thank the organizers for the invitation to speak at the International Colloquium.

\end{document}